\documentclass[11pt, a4paper]{article}
\usepackage{amssymb,amsmath,amsthm}
\usepackage{fancyhdr}
\usepackage[british]{babel}
\usepackage{color}
\usepackage[title]{appendix}
\usepackage{amsmath,amsthm}
\usepackage{enumerate}

\usepackage{hyperref}
\usepackage[margin=2.5cm]{geometry}
\usepackage{lineno}
\usepackage[ruled]{algorithm2e}

\newtheorem{theorem}{Theorem}[section]
\newtheorem{prop}[theorem]{Proposition}
\newtheorem{lemma}[theorem]{Lemma}

\theoremstyle{definition}

\newtheorem{question}{Question}

\newtheorem{defn}[theorem]{Definition}

\newtheorem{claim}[theorem]{Claim}
\newtheorem{fact}[theorem]{Fact}

\newcommand{\btheorem}{\begin{theorem}}
\newcommand{\etheorem}{\end{theorem}}
\newcommand{\bconjecture}{\begin{conjecture}}
\newcommand{\econjecture}{\end{conjecture}}
\newcommand{\bproposition}{\begin{proposition}}
\newcommand{\eproposition}{\end{proposition}}
\newcommand{\bdefinition}{\begin{definition}}
\newcommand{\edefinition}{\end{definition}}
\newcommand{\bcorollary}{\begin{corollary}}
\newcommand{\ecorollary}{\end{corollary}}
\newcommand{\bproof}{\begin{proof}}
\newcommand{\eproof}{\end{proof}}
\newcommand{\bclaim}{\begin{claim}}
\newcommand{\eclaim}{\end{claim}}
\newcommand{\bquestion}{\begin{question}}
\newcommand{\equestion}{\end{question}}
\newcommand{\bfact}{\begin{fact}}
\newcommand{\efact}{\end{fact}}
\newcommand{\bremark}{\begin{remark}}
\newcommand{\eremark}{\end{remark}}
\newcommand{\eexample}{\end{example}}
\newcommand{\bexample}{\begin{example}}

\newcommand{\elemma}{\end{lemma}}
\newcommand{\blemma}{\begin{lemma}}

\title{A robust version of the multipartite Hajnal--Szemer\'edi theorem}

\author{Jie Han\thanks{School of Mathematics and Statistics and Center for Applied Mathematics, Beijing Institute of Technology, Beijing, China. Email: {\tt han.jie@bit.edu.cn}.},
Jie Hu\thanks{Center for Combinatorics and LPMC, Nankai University, Tianjin, China. Email: {\tt hujie@nankai.edu.cn}.},
Donglei Yang\thanks{School of Mathematics, Shandong University, Jinan, China. Email: {\tt dlyang@sdu.edu.cn}.}}

\date{\today}

\begin{document}
\maketitle


\begin{abstract}
  In this note we show the following strengthening of a multipartite version of the Hajnal--Szemer\'edi theorem. For an integer $r\ge 3$ and $\gamma > 0$, there exists a constant $C$ such that if $p\ge Cn^{-2/r}(\log n)^{1/{r \choose 2}}$ and $G$ is a balanced $r$-partite graph with each vertex class of size $n$ and $\delta^\ast(G)\ge (1-1/r+\gamma)n$, then with high probability the random subgraph $G(p)$ of $G$ contains a $K_r$-factor. We also use it to derive corresponding transversal versions.
\end{abstract}

\section{Introduction}
Determining the minimum degree condition for the existence of spanning structures is a central
theme in extremal graph theory, which is called the \emph{Dirac-type} problem because of the cornerstone result of Dirac \cite{Dirac} on the existence of Hamilton cycles.
One of the most frequently studied spanning structures is the \emph{factor}.
Given graphs $H$ and $G$, an \emph{$H$-tiling} of $G$ is a collection of vertex-disjoint copies of $H$ in $G$. An \emph{$H$-factor} of $G$ is an $H$-tiling which covers all vertices of $G$.
The classic theorem of Corr\'adi and Hajnal \cite{Corradi} states the minimum degree condition for the existence of a triangle factor.
More generally, Hajnal and Szemer\'edi \cite{HSz} determined the minimum degree for the existence of a $K_r$-factor.

\begin{theorem}[\cite{Corradi} for the $r=3$ case, \cite{HSz}]\label{thm:HSz}
Let $G$ be an $n$-vertex graph with $n\in r\mathbb{N}$. If $\delta(G)\ge \left(1-1/r\right)n$, then $G$ contains a $K_r$-factor.
\end{theorem}

Let $G$ be an $r$-partite graph with vertex classes $V_1,\ldots,V_r$.
We say that $G$ is \emph{balanced} if $|V_i|=|V_j|$ for any $1\le i < j\le r$.
Write $G[V_i,V_j]$ for the induced bipartite subgraph on vertex classes $V_i$ and $V_j$.
Define $\delta^\ast(G)$ to be $\min_{1\le i<j\le r} \delta(G[V_i,V_j])$.
The following multipartite version of the Hajnal--Szemer\'edi theorem was proved independently and simultaneously by Keevash and Mycroft \cite{MR3290271}, and by Lo and Markstr{\"o}m \cite{Lo}.

\begin{theorem}[\cite{MR3290271,Lo}]\label{thm:multipartite HSz}
For every $\gamma > 0$ and integer $r$, there exists an integer $n_0=n_0(r,\gamma)$ such that if $G$ is a balanced $r$-partite graph with each vertex class of size $n\ge n_0$ and $\delta^\ast(G)\ge (1-1/r+\gamma)n$, then $G$ contains a $K_r$-factor.
\end{theorem}

\subsection{Robustness.}
The study of the \emph{robustness} of graph properties has received considerable attention recently, aiming to enhance classic results in extremal graph theory and probabilistic combinatorics.
We refer the reader to a comprehensive survey of Sudakov \cite{Sudakov-survey} which collect numerous results in this direction.
There are several different measures of robustness mentioned in \cite{Sudakov-survey}, such as random subgraphs, Maker-Breaker games, incompatibility systems and so on.
In this paper the measure of robustness we are interested in is random subgraphs.
Let $0\le p\le 1$ and $G(n,p)$ be the binomial random graph.
Given an $n$-vertex graph $G$, we call $G(p):=G\cap G(n,p)$ the \emph{random subgraph} (or \emph{random sparsification}) of $G$. We say that $G(p)$ has a graph property $\mathcal{P}$ \emph{with high probability}, or \emph{whp} for brevity, if the probability that $G(p)$ has $\mathcal{P}$ tends to 1 as $n$ goes to infinity.
In particular, $G(p)=G(n,p)$ when $G$ is the complete graph $K_n$.

The first result on robustness of graph properties with respect to random subgraphs is a robust version of Dirac's theorem, obtained by Krivelevich, Lee and Sudakov \cite{MR3180741}.
More precisely, they proved that there exists a constant $C$ such that for $p\ge C(\log n)/n$ and an $n$-vertex graph $G$ with $\delta(G)\ge n/2$, whp $G(p)$ contains a Hamilton cycle.
Recently, Allen, B\"ottcher, Corsten, Davies, Jenssen, Morris and Roberts \cite{Allen} proved a robust version of the $r=3$ case of Theorem \ref{thm:HSz}.
Later, Pham, Sah, Sawhney and Simkin \cite{Pham2} derived a robust version of Theorem \ref{thm:HSz} with spread techniques.

\begin{theorem}[\cite{Allen} for the $r=3$ case, \cite{Pham2}]\label{thm:robust HSz}
For any integer $r\ge 3$, there exists a constant $C=C(r)$ such that for any $n\in r\mathbb{N}$ and $p\ge Cn^{-2/r}(\log n)^{1/{r \choose 2}}$ the following holds. If $G$ is an $n$-vertex graph with $\delta(G)\ge \left(1-1/r\right)n$, then whp $G(p)$ contains a $K_r$-factor.
\end{theorem}

For more results on robustness of Dirac-type theorems, we refer the reader to \cite{KKKOP,KMP,Pham2}. In this paper we prove a robust version of Theorem \ref{thm:multipartite HSz}.

\begin{theorem}\label{thm:robust multipartite HSz}
Let $r\ge 3$ be an integer and let $\gamma > 0$. There exists a constant $C=C(r,\gamma)$ such that for any sufficiently large $n$ and $p\ge Cn^{-2/r}(\log n)^{1/{r \choose 2}}$ the following holds. If $G$ is a balanced $r$-partite graph with each vertex class of size $n$ and $\delta^\ast(G)\ge (1-1/r+\gamma)n$, then whp $G(p)$ contains a $K_r$-factor.
\end{theorem}

\subsection{Transversal (robust) versions.}

Recently, there has been much research on the study of \emph{transversal} versions of Dirac-type results \cite{MR4125343,MR3628907,MR4394673,rainbow-bandwidth,CHWW1,CHWW2,MR4287703,Ferber,MR4171383,MR4451911,MR4523452,MR4275007,LLL,MR4451150}.
Given a family of graphs $\mathcal{G}=\{G_1,\ldots,G_m\}$ defined on the same vertex set $V$, an $m$-edge graph $H$ defined on $V$ is \emph{transversal} if $|E(H)\cap E(G_i)|=1$ for each $i\in [m]$. In particular, Cheng, Han, Wang and Wang \cite{CHWW1} proved transversal versions of an asymptotic analogue of the Hajnal--Szemer\'edi theorem and Theorem \ref{thm:multipartite HSz}. It is worth to mention that for general $F$-factor, Montgomery, M\"uyesser and Pehova \cite{MR4451150} derived a transversal version of an asymptotic analogue of the K\"uhn--Osthus theorem \cite{MR2506388}.

\begin{theorem}[\cite{CHWW1}]\label{thm:rainbow HSz}
For every $\gamma > 0$ and integer $r$, there exists $n_0$ such that the following holds for any $n\ge n_0$ and $n\in r\mathbb{N}$. If $G_i$ is an $n$-vertex graph on the vertex set $V$ with $\delta(G_i)\ge (1-1/r+\gamma)n$ for any $i\in[\frac{n}{r}{r\choose2}]$, then $\{G_1,\ldots,G_{\frac{n}{r}{r\choose2}}\}$ contains a transversal $K_r$-factor.
\end{theorem}

\begin{theorem}[\cite{CHWW1}]\label{thm:rainbow multipartite HSz}
For every $\gamma > 0$ and integer $r$, there exists $n_0$ such that the following holds for any $n\ge n_0$. Let $G_i$ be a balanced $r$-partite graph with vertex classes $V_1,\ldots,V_r$ and $|V_j|=n$ for any $i\in[n{r\choose2}]$ and $j\in [r]$. If $\delta^\ast(G_i)\ge (1-1/r+\gamma)n$ for any $i\in[n{r\choose2}]$, then $\{G_1,\ldots,G_{n{r\choose2}}\}$ contains a transversal $K_r$-factor.
\end{theorem}

In this paper, we further consider \emph{transversal robust} versions of Dirac-type theorems, aiming to simultaneously generalize the above results.
Recently, Ferber, Han and Mao \cite{Ferber} studied a transversal robust version of Dirac's theorem and derived Theorem \ref{thm:rainbow robust Dirac} (in fact, they proved a transversal version of resilience result on Hamiltonicity which implies Theorem \ref{thm:rainbow robust Dirac}). Very recently, Anastos and Chakraborti \cite{Anastos} improved this result by removing the error term in the minimum degree condition with spread techniques.

\begin{theorem}[\cite{Ferber}]\label{thm:rainbow robust Dirac}
Let $\gamma > 0$. There exists a constant $C=C(\gamma)$ such that for any sufficiently large $n$ and $p\ge C\log n/n$ the following holds. If $G_i$ is an $n$-vertex graph on the vertex set $V$ with $\delta(G_i)\ge (1/2+\gamma)n$ and $H_i\sim G_i(p)$ for any $i\in[n]$, then whp $\{H_1,\ldots,H_n\}$ contains a transversal Hamilton cycle.
\end{theorem}

The proof of Theorem \ref{thm:rainbow robust Dirac} is based on the local resilience result on Hamiltonicity in random directed graphs by Montgomery \cite{MR4173137}. However, there is no resilience result on $K_r$-factors for $r\ge 3$. Instead of applying the resilience result, we make use of Theorem \ref{thm:robust multipartite HSz} and obtain transversal robust versions of the (multipartite) Hajnal--Szemer\'edi theorem.

\begin{theorem}\label{thm:rainbow robust HSz}
Let $r\ge 3$ be an integer and let $\gamma > 0$. There exists a constant $C=C(r,\gamma)$ such that for any sufficiently large $n\in r\mathbb{N}$ and $p\ge Cn^{-2/r}(\log n)^{1/{r \choose 2}}$ the following holds. If $G_i$ is an $n$-vertex graph on the vertex set $V$ with $\delta(G_i)\ge (1-1/r+\gamma)n$ and $H_i\sim G_i(p)$ for any $i\in[\frac{n}{r}{r\choose2}]$, then whp $\{H_1,\ldots,H_{\frac{n}{r}{r\choose2}}\}$ contains a transversal $K_r$-factor.
\end{theorem}

\begin{theorem}\label{thm:rainbow robust multipartite HSz}
Let $r\ge 3$ be an integer and let $\gamma > 0$. There exists a constant $C=C(r,\gamma)$ such that for any sufficiently large $n$ and $p\ge Cn^{-2/r}(\log n)^{1/{r \choose 2}}$ the following holds. Let $G_i$ be a balanced $r$-partite graph with vertex classes $V_1,\ldots,V_r$ and $|V_j|=n$ for any $i\in[n{r\choose2}]$ and $j\in [r]$. If $\delta^\ast(G_i)\ge (1-1/r+\gamma)n$ and $H_i\sim G_i(p)$ for any $i\in[n{r\choose2}]$, then whp $\{H_1,\ldots,H_{n{r\choose2}}\}$ contains a transversal $K_r$-factor.
\end{theorem}

\section{Notation and preliminaries}
For positive integers $a,b$ with $a<b$, let $[a]=\{1,2,\ldots,a\}$ and $[a,b]=\{a,a+1,\ldots,b\}$.
For constants $x,y,z$, $x=y\pm z$ means that $y-z\le x\le y+z$, and $x\ll y$ means that for any $y> 0$ there
exists $x_0> 0$ such that for any $x< x_0$ the subsequent statement holds.
For a graph $G$, we use $e(G)$ to denote the number of edges in $G$.
For a vertex $v\in V(G)$ and a vertex subset $A$, let $N_A(v)=N_G(v)\cap A$ and $d_A(v)=|N_G(v)\cap A|$.
Define $K_r(G)$ (resp. $K_r(G,v)$) to be the family of all copies of $K_r$ (resp. copies of $K_r$ containing $v$) in $G$.
The \emph{size} of a $K_r$-tiling in $G$ is the number of vertex-disjoint copies of $K_r$ it contains.
The \emph{$r$-clique complex} of a graph $G$ is an $r$-uniform hypergraph whose edges are the $r$-sets of vertices $\{v_1,\ldots,v_r\}$ each inducing a copy of $K_r$ in $G$.
For a vertex subset $A\subseteq V(G)$, let $G[A](p)$ be the random subgraph of the induced subgraph $G[A]$.
Given a set $S$ of size $m$, let $\text{Sym}(S)$ be the set of all $m!$ permutations of $S$.

\subsection{Concentration inequalities}
\begin{lemma}[Chernoff's inequality, \cite{Janson}]\label{lem: chernoff}
Let $X$ be a sum of independent Bernoulli random variables and $\lambda=\mathbb{E}(X)$.
Then for any $0<a<3/2$, we have
$$\mathbb{P}[X\ge (1+a)\lambda]\le e^{-a^2\lambda/3}~~~\text{and}~~~\mathbb{P}[X\le (1-a)\lambda]\le e^{-a^2\lambda/2}.$$
\end{lemma}

\begin{lemma}[Janson's inequality, \cite{Janson}]\label{lem: janson}
Let $G$ be a graph and $\mathcal{F}\subset 2^{E(G)}$ be a collection of subgraphs of $G$ and $p\in[0,1]$. For any $F\in \mathcal{F}$, let $I_F$ be the indicator variable which is 1 if $F$ is present in $G(p)$ and 0 otherwise. Let $X=\sum_{F\in \mathcal{F}}I_F$, $\lambda=\mathbb{E}(X)$ and
$$\bar{\Delta}=\sum_{(F,F')\in \mathcal{F}^2: F\cap F'\neq \emptyset}\mathbb{E}(I_FI_{F'}).$$
Then for any $0<a<1$, we have
$$\mathbb{P}[X\le (1-a)\lambda]\le \exp\left(-\frac{a^2\lambda^2}{2\bar{\Delta}}\right).$$
\end{lemma}

\begin{lemma}[Talagrand-type inequality, \cite{McDiarmid}]\label{lem: permutation ineq}
Let $\{B_1,\ldots,B_k\}$ be a family of finite non-empty sets and $\Omega=\prod_{i}\text{Sym}(B_i)$.
Let $\pi=\{\pi_1,\ldots,\pi_k\}$ be a family of independent permutations, such that for each $i\in [k]$, $\pi_i\in \text{Sym}(B_i)$ is chosen uniformly at random. Let $c$ and $r$ be constants and let $h$ be a nonnegative real-valued function on $\Omega$. Suppose that $h$ satisfies the following conditions for each $\pi \in \Omega$.
\begin{itemize}
  \item If we swap any two elements in any $\pi_i$, the value of $h$ can change  by at most $2c$.
  \item If $h(\pi)=\ell$, then there exists a set $\pi_{proof}\subseteq \pi$ of size at most $r\ell$, such that $h(\pi')\ge \ell$ for any $\pi'\in \Omega$ where $\pi'\supseteq \pi_{proof}$.
\end{itemize}
Let $M$ be the median of $h(\pi)$. Then for each $a\ge 0$, we have
$$P[h(\pi)\le M-a]\le 2\exp\left(-\frac{a^2}{16rc^2M}\right).$$
\end{lemma}
\subsection{Regularity}
Given a graph $G$ and disjoint vertex subsets $X,Y\subseteq V(G)$, the \emph{density} of the pair $(X,Y)$ is defined as $$d(X,Y):=\frac{e(X,Y)}{|X||Y|},$$ where $e(X,Y):=e(G[X,Y])$.
For $\varepsilon>0$, the pair $(X,Y)$ is $\varepsilon$-\emph{regular} if for any $A\subseteq X, B\subseteq Y$ with $|A|\ge \varepsilon |X|, |B|\ge \varepsilon |Y|$, we have $$|d(A,B)-d(X,Y)|<\varepsilon.$$
In addition, $(X,Y)$ is $(\varepsilon,d)$-\emph{regular} if $d(X,Y) \ge d$ for some $d>0$.
An $r$-tuple $(X_1,\ldots,X_r)$ of pairwise disjoint subsets of $V(G)$ is $(\varepsilon,d)$-\emph{regular} if each pair $(X_i,X_j)$ with $1\le i< j\le r$ is $(\varepsilon,d)$-regular.

We say that $(X,Y)$ is $(\varepsilon,d,\delta)$-\emph{super-regular} if $(X,Y)$ is $(\varepsilon,d)$-regular and $d_X(v)\ge \delta |X|$ for any $v\in Y$ and $d_Y(u)\ge \delta |Y|$ for any $u\in X$;
$(X,Y)$ is $(\varepsilon,d)$-\emph{super-regular} if $(X,Y)$ is $(\varepsilon,d,d-\varepsilon)$-super-regular.
The corresponding definitions for regular tuples are analogous.
By the above definitions, we have the following lemma.
\begin{lemma}\label{lem: large degree to super-reg}
Let $0<\varepsilon < 1$ and $(X,Y)$ be a pair of vertex sets such that $d(x,Y)\ge (1-\varepsilon)|Y|$ and $d(y,X)\ge (1-\varepsilon)|X|$ for every $x\in X$ and $y\in Y$. Then $(X,Y)$ is
$(\sqrt{\varepsilon},1-\varepsilon)$-super-regular.
\end{lemma}

\begin{fact}\label{fact: large deg}
Let $(X,Y)$ be an $(\varepsilon,d)$-regular pair, and $B\subseteq Y$ with $|B|\ge \varepsilon |Y|$. Then all but $\varepsilon |X|$ vertices $v\in X$ have $d_B(v)\ge (d-\varepsilon)|B|$.
\end{fact}

\begin{lemma}[Slicing Lemma, \cite{KS}]\label{fact: slicing lem}
Let $(X,Y)$ be an $(\varepsilon,d)$-regular pair, and for some $\eta>\varepsilon$, let $X'\subseteq X, Y'\subseteq Y$ with $|X'|\ge \eta |X|, |Y'|\ge \eta |Y|$. Then $(X',Y')$ is an $\varepsilon'$-regular pair with $\varepsilon'=\max\{\frac{\varepsilon}{\eta},2\varepsilon\}$, and for its density $d'$ we have $d'>d-\varepsilon$.
\end{lemma}

By Fact \ref{fact: large deg} and Lemma \ref{fact: slicing lem}, the following lemma holds.

\begin{lemma}\label{lem: reg to super-reg}
Let $r\in \mathbb{N}$ and $0<\varepsilon <d\le 1$ with $\varepsilon\le\frac{1}{2r}$.
If $(V_1,\ldots,V_r)$ is an $(\varepsilon,d )$-regular tuple with $|V_i|=n$ for every $i\in [r]$, then there are subsets ${V_i'}\subseteq V_i$ with $|{V_i'}|=(1-(r-1)\varepsilon)n$ for every $i\in [r]$ such that the $r$-tuple $({V_1'},\ldots,{V_r'})$ is $(2\varepsilon,d-\varepsilon ,d-r\varepsilon)$-super-regular.
\end{lemma}

\begin{lemma}[Counting Lemma, \cite{KS}]\label{lem: counting lem}
Given $d>\varepsilon>0$, $r,m\in \mathbb{N}$ and an $r$-vertex graph $H$, let $G$ be a graph obtained by replacing every vertex $v_i$ of $H$ with an independent set $V_i$ of size $m$ and every edge $v_iv_j$ of $H$ with an $(\varepsilon,d)$-regular pair $(V_i,V_j)$. If $\varepsilon\le \frac{d^r}{(2+r)2^r}=:d_0$, then there are at least $(d_0m)^r$ copies of $H$ in $G$ so that each $v_i$ is embedded into the set $V_i$.
\end{lemma}

\begin{lemma}[Degree form of the Regularity Lemma, \cite{KS}]\label{lem: reg lem}
Let $G$ be a graph with vertex set $V$.
For every $\varepsilon>0$, there is an $M=M(\varepsilon)$ such that if $d\in[0,1]$ is any real number, then there is a partition $V=V_0\cup V_1\cup \ldots \cup V_k$ and a spanning subgraph $G'\subseteq G$ with the following properties:
\begin{enumerate}
  \item $1/\varepsilon \le k\le M$,
  \item $|V_0|\le \varepsilon |V|$ and $|V_i|=m$ for all $1\le i\le k$ with $m\le \varepsilon |V|$,
  \item $d_{G'}(v)>d_{G}(v)-(d+\varepsilon)|V|$ for all $v\in V$,
  \item $e(G'[V_i])=0$ for all $i\ge 1$,
  \item all pairs $(V_i,V_j)$ $(1\le i<j\le k)$ are $\varepsilon$-regular in $G'$ with density $0$ or at least $d$.
\end{enumerate}
\end{lemma}
\noindent Moreover, $V_0,V_1,\ldots,V_k$ are often referred to as \emph{clusters}.

\begin{defn}[Reduced graph]\label{reduced graph}
Given a graph $G$ with vertex set $V$, a partition $V=V_1\cup \ldots \cup V_k$, and two parameters $\varepsilon,d>0$, the \emph{reduced graph} $R=R(\varepsilon,d)$ of $G$ is defined as follows:
\begin{itemize}
  \item $V(R)=\{V_1,\ldots,V_k\}$,
  \item $V_iV_j\in E(R)$ if and only if $(V_i,V_j)$ is $(\varepsilon,d )$-regular in $G$.
\end{itemize}
\end{defn}

When applying Lemma \ref{lem: reg lem}, we usually begin with a graph $G=(V,E)$ and parameters $\varepsilon,d>0$, and then obtain a partition $V=V_0\cup V_1\cup \ldots \cup V_k$ and a spanning subgraph $G'$ with above-mentioned properties. Then we drop the cluster $V_0$ and study the properties of the reduced graph $R=R(\varepsilon,d)$ of $G'[V_1\cup \ldots \cup V_k]$. By Lemma \ref{lem: reg lem},
\begin{align*}
   \delta(R)\ge \frac{\delta(G)-(d+\varepsilon)|V|-|V_0|}{m} \ge \frac{\delta(G)-(d+2\varepsilon)|V|}{m}.
\end{align*}
In particular, if $\delta(G)\ge c|V|$ for some constant $c$, then $\delta(R)\ge (c-d-2\varepsilon)k$.

\subsection{Random sparsification}
The following result will be used to prove Lemma \ref{lem: rooted lem}, which is stated in \cite{Pham} and a more general form of Riordan and Heckel's result \cite{Heckel,Riordan}.
\begin{lemma}[\cite{Pham}]\label{lem: clique complex}
Let $r\ge 3$. There exists $\varepsilon=\varepsilon(r)>0$ such that for any $p\le n^{-2/r+\varepsilon}$ the following holds. Let $G$ be an $n$-vertex graph and $H$ be the $r$-clique complex of $G$. For some $p_0\sim p^{r\choose 2}$, there is a joint distribution $\lambda$ of a graph $G'$ and an $r$-uniform hypergraph $H'$ defined on the same vertex set $[n]$ such that the following holds. The marginal of $G'$ is the same as a random subgraph $G(p)$ of $G$ and the marginal of $H'$ is the same as a random $r$-uniform subhypergraph $H(p_0)$ of $H$. Furthermore, whp for every hyperedge of $H'$, there is a copy of $K_r$ in $G'$ with the same vertex set.
\end{lemma}

The following lemma generalizes Lemma 8.3 in \cite{Allen}.

\begin{lemma}[Rooted embedding lemma]\label{lem: rooted lem}
For any $0<\mu<\frac{1}{100r}$, there exists a constant $C$ such that the following holds for any sufficiently large $n$ and $p\ge Cn^{-2/r}(\log n)^{1/{r\choose 2}}$.
Let $G$ be an $n$-vertex graph and $v_1,\ldots,v_{\ell}$ be distinct vertices in $G$ with $\ell\le \mu^2n$.
Let $X_1,\ldots,X_t\subseteq V(G)\setminus \{v_1,\ldots,v_{\ell}\}$ be disjoint sets for some $t\in \mathbb{N}$.
If $|K_r(G,v_i)|\ge \mu n^{r-1}$ for each $1\le i\le \ell$, then whp there is a $K_r$-tiling $\{K^1,\ldots,K^{\ell}\}$ in $G(p)$ such that $K^i\in K_r(G,v_i)$ for each $i\in [\ell]$ and $|X_s\cap \bigcup_{i=1}^{\ell}V(K^i)|\le 4r\mu |X_s|+r-2$ for each $s\in [t]$.
\end{lemma}
\begin{proof}[Proof of Lemma \ref{lem: rooted lem}]
Set $C_0=100\mu^{-1}$.
Let $H$ be the $r$-clique complex of $G$ and $p_0={C_0}n^{1-r}\log n$.
For each $i\in [\ell]$, define $E(H,v_i)$ to be the set of hyperedges in $H$ containing $v_i$. Then $|E(H,v_i)|=|K_r(G,v_i)|\ge \mu n^{r-1}$.
By Lemma \ref{lem: clique complex}, it suffices to prove that whp there is a matching $\{K^1,\ldots,K^{\ell}\}$ in $H(p_0)$ such that $K^i\in E(H,v_i)$ for each $i\in [\ell]$ and $|X_s\cap \bigcup_{i=1}^{\ell}V(K^i)|\le 4r\mu |X_s|+r-2$ for each $s\in [t]$.

We will reveal hyperedges of $H$ step by step and choose hyperedges $K^1,\ldots,K^{\ell}$ one at a time.
A hyperedge $e\in E(H)$ is \emph{active} if it has not been revealed yet.
Given $s\in[t]$ and $i\in [\ell]$, $X_s$ is \emph{saturated} at time $i$ if $|X_s\cap \bigcup_{j=1}^{i-1}V(K^j)|\ge 4r\mu |X_s|$. Let $Y_i$ be the union of the saturated sets $X_s$ at time $i$.
Since $X_s$'s are disjoint, we have
$|Y_i|\le \sum_{s\in [t]}\frac{|X_s\cap \bigcup_{j=1}^{i-1}V(K^j)|}{4r\mu}\le  \frac{|\bigcup_{j=1}^{\ell}V(K^j)|}{4r\mu} =\frac{r\ell}{4r\mu}=\frac{\ell}{4\mu}$.

For each $i$, we choose an arbitrary set $E'(H,v_i)\subseteq E(H,v_i)$ of size exactly $\mu n^{r-1}$. In step $i$, we reveal all \emph{candidate} hyperedges in $E'(H,v_i)$, which are active but not containing any vertex in $\bigcup_{j=1}^{i-1}V(K^j)\cup \{v_{i+1},\ldots,v_{\ell}\}\cup Y_i$.
We claim that in the beginning of step $i$ whp the number of candidate hyperedges in $E'(H,v_i)$ is at least $\frac{\mu}{2} n^{r-1}$.
Indeed, the number of hyperedges that have been revealed before step $i$ is at most $2\mu n^{r-1}p_0(i-1)\le2\mu n^{r-1}p_0\ell\le 2\mu^3C_0n\log n\le \frac{\mu}{8}n^{r-1}$ with probability at least $1-(i-1) n^{-2\mu C_0/3}$ by the Chernoff's inequality and a union bound.
Moreover, the number of edges both containing $v_i$ and some vertex in $\bigcup_{j=1}^{i-1}V(K^j)\cup \{v_{i+1},\ldots,v_{\ell}\}\cup Y_i$ is at most $\left(r\ell+\frac{\ell}{4\mu}\right)\cdot n^{r-2}\le \left(r\mu^2n+\frac{\mu n}{4}\right)\cdot n^{r-2}\le \frac{3\mu}{8}n^{r-1}$.

Now we reveal all candidate hyperedges in $E'(H,v_i)$.
Since these candidates are pairwise independent, by the Chernoff's inequality, there exist at least $\frac{\mu}{4} n^{r-1}p_0=\frac{1}{4}\mu C_0\log n~(\ge 1)$ hyperedges in $E(H(p_0))\cap E'(H,v_i)$ with probability at least $1-n^{-\mu C_0/16}$.
We arbitrarily take one of them as $K^i$.
Then move on to the next step $i+1$.
Finally, by a union bound, whp there is a desired matching $\{K^1,\ldots,K^{\ell}\}$ in $H(p_0)$ with $K^i\in E(H,v_i)$.
\end{proof}

\begin{prop}\label{prop: standard computation}
Let $r\ge 3$.
For any $\alpha>0$, there exists a constant $C$ such that the following holds
for any $n\in \mathbb{N}$ and $p\ge Cn^{-2/r}(\log n)^{1/{r \choose 2}}$.
Let $G$ be an $n$-vertex graph and $\mathcal{F}$ be a family of copies of $K_r$ in $G$.
We denote by $\mathcal{F}(p)$ the family of copies of $K_r$ in $\mathcal{F}$
that are present in $G(p)$.
If $|\mathcal{F}|= \alpha n^{r}$, then with probability at least $1-n^{-\alpha C^{r \choose 2}n/16}$, $|\mathcal{F}(p)|\ge \frac{1}{2}\alpha C^{r\choose 2}n\log n$.
\end{prop}

\begin{proof}[Proof of Proposition \ref{prop: standard computation}]
Assume that $p=Cn^{-2/r}(\log n)^{1/{r \choose 2}}$.
For each $F\in \mathcal{F}$, let $I_F$ be the indicator random variable for the event $F\in \mathcal{F}(p)$ and $X=\sum_{F\in \mathcal{F}}I_F$. So $|\mathcal{F}(p)|=X$. We have  $\mathbb{E}(X)=|\mathcal{F}|\cdot p^{r \choose 2}= \alpha C^{r\choose 2}n\log n$. To calculate $\bar{\Delta}=\sum_{(F,F')\in \mathcal{F}^2: F\cap F'\neq \emptyset}\mathbb{E}(I_FI_{F'})$, first note that for $2\le s\le r-1$, we have
$$n^{s-1}p^{s \choose 2}\ge n^{s-1-\frac{2}{r}{s \choose 2}}=n^{(s-1)(1-\frac{s}{r})}\ge n^{1-2/r}.$$
Thus,  we have
\begin{align*}
\bar{\Delta}-\mathbb{E}(X)
&\le \sum_{s=2}^{r-1}|\mathcal{F}|\cdot{r\choose s} n^{r-s}p^{2 {r\choose 2}-{s\choose 2}}
\\& \le r^{r-1}\sum_{s=2}^{r-1} \alpha n^{2r-s}p^{2 {r\choose 2}-{s\choose 2}}
\\& =  r^{r-1}\alpha \sum_{s=2}^{r-1} n^{2r-1}p^{2 {r\choose 2}}\cdot (n^{s-1}p^{s\choose 2})^{-1}
\\& \le  r^{r} \alpha C^{r^2-r}(\log n)^2 n^{2/r}
\\& = o(\mathbb{E}(X)).
\end{align*}
By Janson's inequality, with probability at least $1-n^{-\alpha C^{r \choose 2}n/16}$, $X\ge \mathbb{E}(X)/2$, that is, $|\mathcal{F}(p)|\ge \frac{1}{2}\alpha C^{r\choose 2}n\log n$.
\end{proof}

The next lemma generalizes Lemma 8.1 (ii) in \cite{Allen}.

\begin{lemma}\label{lem: vtx-disjoint Kr}
For any $\mu>0$, there exists a constant $C$ such that the following holds for any $n,n_1,n_2\in \mathbb{N}$ with $n_1\ge n_2$ and $p\ge Cn^{-2/r}(\log n)^{1/{r \choose 2}}$. Let $G=G(V_1,\ldots,V_r,E)$ be an $n$-vertex $r$-partite graph with $|V_i|\ge n_1$ for each $i\in[r]$. If for any $X_i\subseteq V_i$ with $|X_i|\ge n_2$, $G[X_1,\ldots,X_r]$ contains at least $\mu {n^r}$ copies of $K_r$, then whp $G(p)$ contains at least $n_1-n_2$ vertex-disjoint copies of $K_r$.
\end{lemma}

\begin{proof}[Proof of Lemma \ref{lem: vtx-disjoint Kr}]
It suffices to prove that whp there exists a copy of $K_r$, even after deleting any set $U$ of size $r(n_1-n_2-1)$ with $|U\cap V_i|=n_1-n_2-1$ for each $i\in[r]$. Fix one such $U$. Let $X_i=V_i\setminus U$ for each $i\in [r]$. Then $|X_i|\ge n_2+1$ for all $i$ and so $G[X_1,\ldots,X_r]$ contains at least $\mu n^r$ copies of $K_r$. We arbitrarily choose a family of $\mu n^{r}$ copies of $K_r$ from $G[X_1,\ldots,X_r]$, denoted by $\mathcal{F}$. We apply Proposition \ref{prop: standard computation} on $\mathcal{F}$ with $\alpha=\mu$ and conclude that with probability at least $1-n^{-\mu C^{r \choose 2}n/16}$, $|\mathcal{F}(p)|\ge \frac{1}{2}\mu C^{r\choose 2}n\log n>0$.

By a union bound, the probability that there exists a set $U$ of size $r(n_1-n_2-1)$ and $|U\cap V_i|=n_1-n_2-1$ for each $i\in[r]$ such that $|\mathcal{F}(p)|=0$ is at most
$$2^n\cdot n^{-\mu C^{r \choose 2}n/16}=o(1).$$
Hence whp there exists a copy of $K_r$ in $G(p)$, even after deleting any set $U$ of size $r(n_1-n_2-1)$ with $|U\cap V_i|=n_1-n_2-1$ for each $i\in[r]$.
Therefore, whp $G(p)$ contains at least $n_1-n_2$ vertex-disjoint copies of $K_r$.
\end{proof}

The following lemma can be obtained by Theorem 1.9 in \cite{Pham2}, Theorem 1.6 in \cite{Frankston} and Lemma \ref{lem: clique complex}, following the proof strategy of Theorem 1.6 in \cite{Pham2}, and its proof will be included in Appendix \ref{appendix} for completeness.

\begin{lemma}\label{lem: Pham}
Let $r\ge 3$ and $d>0$. There exist constants $\varepsilon,C>0$ such that the following holds for any sufficiently large $n$ and $p\ge Cn^{-2/r}(\log n)^{1/{r \choose 2}}$. Let $G=(V_1,V_2,\ldots, V_r,E)$ be a balanced $r$-partite graph with each vertex class of size $n$. If $(V_i,V_j)$ is $(\varepsilon,d_{ij})$-super regular with $d_{ij}\ge d$ for all $i\neq j$, then whp $G(p)$ contains a $K_r$-factor.
\end{lemma}

The next lemma gives a relation between weight functions on the vertex set of a graph $G$ and the family of all copies of $K_r$ in $G$, and the proof idea follows Corollary 7.3 in \cite{Allen}.

\begin{lemma}\label{lem: integer weight}
Let $r\ge 2$ be an integer and let $\gamma > 0$. There exists an integer $n_0=n_0(r,\gamma)$ such that the following holds for any $n\ge n_0$.
Let $G=(V_1,V_2,\ldots,V_r,E)$ be a balanced $r$-partite graph with each vertex class of size $n$ and $\delta^\ast(G)\ge (1-1/r+\gamma)n$.
Let $\lambda: V(G)\rightarrow \mathbb{N}_+$ be a weight function such that
$\sum_{v\in V_i}\lambda(v)=\sum_{u\in V_j}\lambda(u)$ for any $1\le i<j\le r$ and $\lambda(v)=\left(1\pm\frac{\gamma}{4}\right)\frac{\sum_{u\in V_i}\lambda(u)}{n}$ for any $v\in V_i$ and $i\in [r]$.
Then there exists a weight function $\omega: K_r(G)\rightarrow \mathbb{N}_0$ such that $\sum_{K\in K_r(G,v)}\omega(K)=\lambda(v)$ for any $v\in V(G)$.
\end{lemma}

\begin{proof}[Proof of Lemma \ref{lem: integer weight}]
Given an integer $r\ge 2$ and $\gamma>0$, we choose $$\frac{1}{n}\ll \gamma, \frac{1}{r}.$$
We first define an auxiliary graph $H$ by replacing every vertex $v\in V(G)$ with an independent set $I_v$ of size $\lambda(v)$ and every edge $uv\in E(G)$ is replaced by a complete bipartite graph with two parts $I_u,I_v$.
Since $\sum_{v\in V_i}\lambda(v)=\sum_{u\in V_j}\lambda(u)$ for any $1\le i<j\le r$, $H$ is a balanced $r$-partite graph with each vertex class of size $\sum_{v\in V_i}\lambda(v)$.
Let $N=\sum_{v\in V_i}\lambda(v)$.
Then $\delta^\ast(H)\ge (1-1/r+\gamma)n \cdot\left(1-\frac{\gamma}{4}\right)\frac{\sum_{u\in V_i}\lambda(u)}{n}= (1-1/r+\gamma)\left(1-\frac{\gamma}{4}\right)N> (1-1/r+\gamma/2)N$.
By Theorem \ref{thm:multipartite HSz}, there exists a weight function $\omega_H: K_r(H)\rightarrow \{0,1\}$ such that $\sum_{K\in K_r(H,v)}\omega(K)=1$ for any $v\in V(H)$. Now we define $\omega: K_r(G)\rightarrow \mathbb{N}_0$ by $\omega(K)=\sum_{K'\in K_r(H[K])}\omega_H(K')$, where $H[K]$ is the subgraph of $H$ induced by $\cup_{v\in V(K)}I_v$. Hence for every $v\in V(G)$,
\begin{align*}
\sum_{K\in K_r(G,v)}\omega(K)
&=\sum_{K\in K_r(G,v)}\sum_{K'\in K_r(H[K])}\omega_H(K')
\\& =\sum_{v'\in I_v}\sum_{K'\in K_r(H,v')}\omega_H(K')
\\& =\sum_{v'\in I_v}1=|I_v|=\lambda(v).
\end{align*}
\end{proof}

\section{Proof of Theorem \ref{thm:robust multipartite HSz}}
We follow the proof of Theorem 1.2 in \cite{Allen} for the case no large sparse set. Given an integer $r\ge 2$ and $\gamma>0$, we choose $$\frac{1}{n}\ll\frac{1}{C}\ll\varepsilon\ll d,\alpha\ll \gamma, \frac{1}{r}.$$
Recall that $G$ is a balanced $r$-partite graph with each vertex class of size $n$.
Let $V_1,\ldots,V_r$ be the $r$ vertex classes of $G$.
By applying Lemma \ref{lem: reg lem} on $G$ with constants $\varepsilon$ and $d$, we obtain a partition $V_0\cup\bigcup_{i\in[r],j\in[k_i]}V_{ij}$ such that  $V_{ij} \subseteq V_i$ and $\frac{1}{\varepsilon}\le \sum_{i=1}^r k_i\le M$, which refines the partition $V_1\cup\ldots\cup V_r$.
We may assume that $k_1=\cdots=k_r=:k$ by placing clusters into $V_0$ and updating $V_0$ if necessary.
Then $|V_0|\le r\cdot \varepsilon |G|=\varepsilon r^2 n$ and $|V_{ij}|\ge \frac{n-|V_0|}{k}\ge \frac{n-\varepsilon r^2n}{k}$ for all $i\in[r],j\in[k]$.
Let $R:=R(\varepsilon,d)$ be the reduced graph for this partition.
Then $R$ is a balanced $r$-partite graph with each vertex class of size $k$.
By the choice $\varepsilon\ll d\ll \gamma$ and the fact that $\delta^\ast(G)\ge (1-1/r+\gamma)n$, we have $\delta^\ast(R)\ge (1-1/r+\gamma/2)k$. By Theorem \ref{thm:multipartite HSz}, $R$ contains a $K_r$-factor $\{K^1,K^2,\ldots,K^k\}$.
Without loss of generality, let $V(K^j)=\{V_{ij}:1\le i\le r\}$ for each $j\in [k]$.
Then the $r$-tuple $(V_{1j},V_{2j},\ldots,V_{rj})$ is $(\varepsilon,d)$-regular for each $j\in [k]$.
By Lemma \ref{lem: reg to super-reg}, there are subsets $V_{ij}'\subseteq V_{ij}$ with $|V_{ij}'|=(1-(r-1)\varepsilon)|V_{ij}|$ for all $i\in [r]$ so that the $r$-tuple $(V_{1j}',V_{2j}',\ldots,V_{rj}')$ is $(2\varepsilon,d-\varepsilon,d-r\varepsilon)$-super-regular for each $j\in [k]$.
Let $B=V_0\cup \bigcup_{i\in[r],j\in[k]}(V_{ij}\setminus V_{ij}')$.
Then $|B|\le r\cdot\varepsilon |G|+\sum_{i\in[r],j\in[k]}(r-1)\varepsilon|V_{ij}|\le 2\varepsilon r^2 n$.

Recall that $p\ge Cn^{-2/r}(\log n)^{1/{r \choose 2}}$.
Let $p'$ be such that $1-p=(1-p')^3$.
We will sprinkle edges of $G(p)$ in three rounds.
We first find a $K_r$-tiling covering vertices in $B$ by Lemma \ref{lem: rooted lem}.
Then we remove some vertex-disjoint copies of $K_r$ to balance the super-regular $r$-tuples by Lemma \ref{lem: integer weight}.
Finally we apply Lemma \ref{lem: Pham} on the balanced super-regular $r$-tuples.

\vbox{}
\noindent \textbf{Cover vertices in $B$.}

Note that $V(G)\setminus B=\cup_{i\in[r],j\in[k]}V_{ij}'$.
By the choice $\alpha\ll \gamma$ and choosing each vertex of $V(G)\setminus B$ independently with probability $1/2$, we can find a set $W \subset V(G)\setminus B$ such that
\begin{enumerate}[(1)]
  \item $|W\cap V_{ij}'|=\left(\frac{1}{2}\pm \alpha\right)\frac{n}{k}$ for any $i\in[r],j\in[k]$;
  \item\label{2} For any $i\in [r]$, $d_G(v,W\cap V_i)\ge \left(1-\frac{1}{r}+\frac{\gamma}{4}\right)|W\cap V_i|$ for each $v\in V(G)\setminus V_i$;
  \item\label{3} $d_G(v,V_{ij}'\cap W)=\left(\frac{1}{2}\pm\frac{1}{4}\right)d_G(v,V_{ij}')$ for each $v\in V(G)$ with $d_G(v,V_{ij}')\ge \varepsilon |V_{ij}'|$.
\end{enumerate}
When covering vertices in $B$, we only use vertices in $W$ to preserve the super-regularity of the $r$-tuples $(V_{1j}',V_{2j}',\ldots,V_{rj}')$.
\begin{claim}\label{claim: cover vertices in B}
Whp there is a $K_r$-tiling $\mathcal{K}_1\subseteq G[B\cup W]$ in $G(p')$ such that $B\subseteq V(\mathcal{K}_1)$ and $|V(\mathcal{K}_1)\cap V_{ij}'|\le 8r^2\sqrt{\varepsilon}|V_{ij}'|$ for all $i\in[r],j\in[k]$.
\end{claim}
\begin{proof}[Proof of Claim \ref{claim: cover vertices in B}]
Since $|B|\le 2 \varepsilon r^2 n$ and $|W|\ge rn/4$, we have $|B|\le 8\varepsilon r |B\cup W|$.
For each $v\in B$, we denote by $\mathcal{F}(v)$ the family of copies of $K_r$ containing $v$ and at most one vertex of $W\cap V_i$ for each $i\in [r]$.
By (\ref{2}) above, given $i\in [r]$, any $r-1$ vertices of $V(G)\setminus V_i$ have at least $(r-1)\left(1-\frac{1}{r}+\frac{\gamma}{4}\right)|W\cap V_i|-(r-2)|W\cap V_i|\ge \left(\frac{1}{r}+\frac{\gamma(r-1)}{4}\right)\frac{n}{4}$ common neighbours in $W\cap V_i$.
Hence $|\mathcal{F}(v)|\ge \left(\frac{n}{4r}\right)^{r-1}$ for each $v\in B$.
By applying Lemma \ref{lem: rooted lem} on $B$ with $\mu=2\sqrt{2\varepsilon r}$ and $X_s=V_{ij}'$, the claim follows.
\end{proof}

\noindent \textbf{Balance super-regular $r$-tuples.}

Let $V_{ij}''=V_{ij}'\setminus V(\mathcal{K}_1)$ for all $i\in[r],j\in[k]$.
Now for each $j\in [k]$, the super-regular $r$-tuple $(V_{1j}'',V_{2j}'',\ldots,V_{rj}'')$ may not be balanced. We will balance them by the following claim.
\begin{claim}\label{claim: balnace super-regular tuple}
Whp there is a $K_r$-tiling $\mathcal{K}_2\subseteq G[W\setminus V(\mathcal{K}_1)]$ in $G(p')$ such that $|V_{ij}'' \setminus V(\mathcal{K}_2)|=\lfloor\frac{9n}{10k}\rfloor$ for all $i\in[r],j\in[k]$.
\end{claim}
\begin{proof}[Proof of Claim \ref{claim: balnace super-regular tuple}]
We define a weight function $\lambda: V(R)\rightarrow \mathbb{N}_+$ by $\lambda(V_{ij})=|V_{ij}''|-\lfloor\frac{9n}{10k}\rfloor$.
Then by Claim \ref{claim: cover vertices in B}, $\left(\frac{1}{10}-10r^2\sqrt{\varepsilon}\right)\frac{n}{k} \le \lambda(V_{ij})\le \lceil\frac{n}{10k}\rceil$
and $\sum_{j=1}^k\lambda(V_{ij})=n-|\mathcal{K}_1|- k\lfloor\frac{9n}{10k}\rfloor$ for each $i\in[r]$.
Recall that $\delta^\ast(R)\ge (1-1/r+\gamma/2)k$.
By applying Lemma \ref{lem: integer weight} on $R$ and $\lambda$, there exists a weight function $\omega: K_r(R)\rightarrow \mathbb{N}_0$ such that $\sum_{K\in K_r(R,V_{ij})}\omega(K)=\lambda(V_{ij})$ for any $V_{ij}\in V(R)$.

For every $K\in K_r(R)$, we shall find a family $\mathcal{F}_K$ of $\omega(K)$ vertex-disjoint copies of $K_r$ in $G[W\cap \bigcup_{V_{ij}\in V(K)}V_{ij}''](p')$ such that all $\mathcal{F}_K$'s are pairwise vertex-disjoint, and then let $\mathcal{K}_2=\cup_{K\in K_r(R)}\mathcal{F}_K$.
We shall iteratively build $\mathcal{F}_K$ by going through every $K\in K_r(R)$ in order.
Observe that at the end of each step, the set $L_{ij}$ of uncovered vertices in $V_{ij}''\cap W$ has order at least $|V_{ij}''\cap W|-\lambda(V_{ij})=|V_{ij}'\cap W|-|V(\mathcal{K}_1)\cap V_{ij}'|-\lambda(V_{ij})\ge \frac{3n}{10k}$ for all $i\in[r],j\in[k]$.
We claim that we can find such families $\mathcal{F}_K$ for every $K\in K_r(R)$.
In fact, for each fixed $K$, without loss of generality, let $V(K)=\{V_{11},V_{21},\ldots,V_{r1}\}$.
By Lemma \ref{fact: slicing lem}, for any $X_{i}\subseteq V_{i1}''\cap W$, $X_{j}\subseteq V_{j1}''\cap W$ with $|X_{i}|, |X_{j}|\ge d\cdot\frac{n}{k}$, we have $(X_{i},X_{j})$ is $\left(\frac{\varepsilon}{d}, d-\varepsilon\right)$-regular for all $i\neq j\in [r]$.
Then by Lemma \ref{lem: counting lem}, $|K_r(G[X_{1},\ldots,X_{r}])|\ge\left(\frac{(d-\varepsilon)^r}{(2+r)2^r}\cdot\frac{dn}{k}\right)^r$.
Since $|L_{i1}|\ge \frac{3n}{10k}$ for any $i\in [r]$,
by applying Lemma \ref{lem: vtx-disjoint Kr} with $n_1=\frac{3n}{10k}$ and $n_2=\frac{dn}{k}$, whp $G[\cup_{i\in [r]}L_{i1}](p')$ contains at least $\left(\frac{3}{10}-d\right)\frac{n}{k}>\frac{n}{5k}=2\cdot\frac{n}{10k}$ vertex-disjoint copies of $K_r$.
Recall that $\omega(K)\le \lambda(V_{i1})\le \lceil\frac{n}{10k}\rceil$ for any $i\in [r]$.
Hence we can find the desired family $\mathcal{F}_K$ for $K$.
Recall that $\mathcal{K}_2=\cup_{K\in K_r(R)}\mathcal{F}_K$.
Then $|V_{ij}'' \setminus V(\mathcal{K}_2)|=\lfloor\frac{9n}{10k}\rfloor$ for all $i\in[r],j\in[k]$.
\end{proof}

\noindent \textbf{Apply Lemma \ref{lem: Pham} on the balanced super-regular $r$-tuples.}

Let $U_{ij}=V_{ij}'' \setminus V(\mathcal{K}_2)$ for all $i\in [r], j\in [k]$.
Then $|U_{ij}|=\lfloor\frac{9n}{10k}\rfloor$ for all $i\in[r],j\in[k]$.
Recall that $(V_{1j}',V_{2j}',\ldots,V_{rj}')$ is $(2\varepsilon,d-\varepsilon,d-r\varepsilon)$-super-regular for every $j\in [k]$.
By (\ref{3}) above, for all $i_1,i_2\in[r]$, $j\in [k]$ and $v\in U_{i_1j}$, we have $d_G(v,U_{i_2j})\ge d_G(v,V_{i_2j}'\setminus W)\ge \frac{1}{4}d_G(v,V_{i_2j}')\ge \frac{d-r\varepsilon}{4}|V_{i_2j}'|\ge \frac{d}{8}|V_{i_2j}'|\ge \frac{d}{8}|U_{i_2j}|$.
By Lemma \ref{fact: slicing lem},  $(U_{1j},U_{2j},\ldots,U_{rj})$ is $(\frac{20\varepsilon}{9},\frac{d}{2},\frac{d}{8})$-super-regular for each $j\in [k]$, which is also $(\frac{20\varepsilon}{9},\frac{d}{8})$-super-regular.
By applying Lemma \ref{lem: Pham}, whp $G[U_{1j},U_{2j},\ldots,U_{rj}](p')$ contains a $K_r$-factor for every $j\in [k]$, and they form a $K_r$-factor of $G[\cup_{i\in [r],j\in [k]}U_{ij}]$, denoted by $\mathcal{K}_3$.

In summary, $\mathcal{K}_1\cup\mathcal{K}_2\cup\mathcal{K}_3$ is a $K_r$-factor in $G(p)$ as desired.

\section{Proofs of Theorems \ref{thm:rainbow robust HSz} and \ref{thm:rainbow robust multipartite HSz}}


\begin{proof}[Proof of Theorem \ref{thm:rainbow robust HSz} from Theorem \ref{thm:rainbow robust multipartite HSz}]
We first claim that there is a balanced partition $V=V_1\cup\ldots\cup V_r$ such that $d_{G_i}(v,V_j)\ge \left(1-\frac{1}{r}+\frac{\gamma}{2}\right)\frac{n}{r}$ for any $v\in V$, $i\in [\frac{n}{r}{r\choose2}]$ and $j\in[r]$. Indeed, choose a balanced partition $V=V_1\cup\ldots\cup V_r$ uniformly at random. For fixed $v\in V$, $i\in [\frac{n}{r}{r\choose2}]$ and $j\in[r]$, $\mathbb{E}[d_{G_i}(v,V_j)]=\frac{d_{G_i}(v)}{r}$.
By Chernoff's inequality, $\mathbb{P}[d_{G_i}(v,V_j)\le \left(1-\frac{\gamma}{2}\right)\frac{d_{G_i}(v)}{r}]\le \exp\left(-\frac{\gamma^2}{8}\frac{d_{G_i}(v)}{r}\right)<\exp\left(-\frac{\gamma^2}{16r}n\right)$.
By a union bound, 
whp $d_{G_i}(v,V_j)\ge \left(1-\frac{1}{r}+\frac{\gamma}{2}\right)\frac{n}{r}$ for any $v\in V$, $i\in [\frac{n}{r}{r\choose2}]$ and $j\in[r]$.

Recall that $H_i\sim G_i(p)$ for any $i\in[\frac{n}{r}{r\choose2}]$. Let $G_i'=G_i[V_1,\ldots,V_r]$ for any $i\in [\frac{n}{r}{r\choose2}]$.
Then $\delta^\ast(G_i')\ge \left(1-\frac{1}{r}+\frac{\gamma}{2}\right)\frac{n}{r}$ for any $i\in [\frac{n}{r}{r\choose2}]$.
By Theorem \ref{thm:rainbow robust multipartite HSz}, if $H_i'\sim G_i'(p)$ for any $i\in[\frac{n}{r}{r\choose2}]$, then whp $\{H_1',\ldots,H_{\frac{n}{r}{r\choose2}}'\}$ contains a transversal $K_r$-factor. So does $\{H_1,\ldots,H_{\frac{n}{r}{r\choose2}}\}$ as $H_i'\subseteq H_i$ for every $i$.
\end{proof}

Before the proof of Theorem \ref{thm:rainbow robust multipartite HSz} starts, we first define an auxiliary $r$-partite graph and the idea comes from \cite{Ferber}.
Recall that $G_i$ is a balanced $r$-partite graph with vertex classes $V_1,\ldots,V_r$ and $|V_j|=n$ for any $i\in[n{r\choose2}]$ and $j\in [r]$.
Let $L:=\{(i,j):i,j\in [r]~\text{and}~i<j\}$.
We divide $\{G_i:i\in[n{r\choose2}]\}$ into ${r\choose 2}$ subfamilies $\{\mathcal{F}_{i,j}\}_{(i,j)\in L}$ with $\mathcal{F}_{i,j}=\{G_{c_{i,j}+1},\ldots,G_{c_{i,j}+n}\}$, where $c_{i,j}=\sum_{(s,t)\prec (i,j)}|\mathcal{F}_{s,t}|$ and $\prec$ is the lexicographic ordering on $L$.
Without loss of generality, let $V_i=[(i-1)n+1,in]$ for each $i\in [r]$.
Let $\pi_i$ be a uniformly random permutation of $V_i$ for each $i\in [r]$.
Now we define auxiliary bipartite graphs $B_{\pi_i}[V_i,V_j]$ for any $1\le i<j\le r$: For any $s\in V_i$ and $t\in V_j$, 
we connect $s$ and $t$ if and only if $st\in E(G_{c_{i,j}+\pi_i(s)-(i-1)n})$.
Note that $G_{c_{i,j}+\pi_i(s)-(i-1)n}\in \mathcal{F}_{i,j}$.
Let $\pi:=\{\pi_1,\ldots,\pi_r\}$ and $B_{\pi}:=\mathop{\cup}\limits_{1\le i<j\le r}B_{\pi_i}[V_i,V_j]$.

\begin{proof}[Proof of Theorem \ref{thm:rainbow robust multipartite HSz}]
Let $B_{\pi}$ be the $r$-partite graph defined as above.
Note that a $K_r$-factor in $B_{\pi}$ (resp. $B_{\pi}(p)$) corresponds to a transversal $K_r$-factor in $\{G_1,\ldots,G_{n{r\choose2}}\}$ (resp. $\{H_1,\ldots,H_{n{r\choose2}}\}$).
Hence, if $\delta^\ast(B_{\pi})\ge \left(1-\frac{1}{r}+\frac{\gamma}{2}\right)n$, then by Theorem \ref{thm:robust multipartite HSz}, whp $B_{\pi}(p)$ contains a $K_r$-factor and Theorem \ref{thm:rainbow robust multipartite HSz} follows. The rest of the proof is to show that whp $\delta^\ast(B_{\pi})\ge \left(1-\frac{1}{r}+\frac{\gamma}{2}\right)n$.

\begin{claim}[\cite{Ferber}, Lemma 2.14]\label{claim: median}
Let $0<\alpha<\frac{1}{2}$. Let $\pi_1$ be a uniformly random permutation on $V_1$ and $\lambda_t=\mathbb{E}[d_{B_{\pi_1}[V_1,V_2]}(t)]$ for every $t\in V_2$. If $\delta(G_i[V_1,V_2])\ge 200/\alpha^2$ for all $i\in [n]$, then we have
$$M_t:=M(d_{B_{\pi_1}[V_1,V_2]}(t))= (1\pm \alpha)\lambda_t.$$
\end{claim}

\begin{claim}\label{claim: permutation}
For any $1\le i<j\le r$ and a uniformly random permutation $\pi_i$ of $V_i$, whp $\delta(B_{\pi_i}[V_i,V_j])\ge \left(1-\frac{1}{r}+\frac{\gamma}{2}\right)n$.
Furthermore, whp $\delta^\ast(B_{\pi})\ge \left(1-\frac{1}{r}+\frac{\gamma}{2}\right)n$.
\end{claim}
\begin{proof}[Proof of Claim \ref{claim: permutation}]
By symmetry, we only show that for a uniformly random permutation $\pi_1$ of $V_1$, whp $\delta(B_{\pi_1}[V_1,V_2])\ge \left(1-\frac{1}{r}+\frac{\gamma}{2}\right)n$, and then apply a union bound on all $r\choose 2$ bipartite subgraphs.
Since $\delta (G_{\pi_1(s)}[V_1,V_2])\ge \left(1-\frac{1}{r}+\gamma\right)n$, we have $d_{B_{\pi_1}[V_1,V_2]}(s)\ge \left(1-\frac{1}{r}+\gamma\right)n$ for any $s\in V_1$.
For a fixed vertex $t\in V_2$, let $\lambda_t=\mathbb{E}[d_{B_{\pi_1}[V_1,V_2]}(t)]$. Then $\lambda_t\ge \left(1-\frac{1}{r}+\gamma\right)n$. We will show that whp $d_{B_{\pi_1}[V_1,V_2]}(t)\ge \left(1-\frac{1}{r}+\frac{\gamma}{2}\right)n$ by applying Lemma \ref{lem: permutation ineq} with $k=1$, $B_1=V_1$ and $h(\pi_1)=d_{B_{\pi_1}[V_1,V_2]}(t)$.
Note that if we swap any two elements in $\pi_1$, the value of $h(\pi_1)$ can change by at most $2$. Additionally, if $h(\pi_1)=d_{B_{\pi_1}[V_1,V_2]}(t)=\ell$, then we choose $\pi_{proof}$ as the $\ell$ indices reflected in $N_{B_{\pi_1}[V_1,V_2]}(t)$.
Hence, $h(\pi_1)$ satisfies the conditions of Lemma \ref{lem: permutation ineq} with $c=1$ and $r=1$.

Now let $\alpha=\gamma/100$.
Then by Claim \ref{claim: median}, the medium $M_t$ of $d_{B_{\pi_1}[V_1,V_2]}(t)$ is at least $(1-\alpha)\lambda_t$.
Hence $\left(1-\frac{1}{r}+\frac{\gamma}{2}\right)n\le \left(1-\frac{\gamma}{2}\right)\lambda_t \le \left(1-\frac{\gamma}{2}\right)\frac{M_t}{1-\alpha}\le \left(1-\frac{\gamma}{4}\right)M_t.$ Thus by Lemma ~\ref{lem: permutation ineq}, we have
\begin{align*}
       & \mathbb{P}\left[h(\pi_1)\le \left(1-\frac{1}{r}+\frac{\gamma}{2}\right)n\right]
\\ \le & \mathbb{P}\left[h(\pi_1)\le \left(1-\frac{\gamma}{4}\right)M_t\right]
\\ \le & 2\exp\left\{-\frac{(\gamma M_t/4)^2}{16M_t}\right\}
\\ \le & 2\exp\left\{-\left(\frac{\gamma}{32}\right)^2n\right\}.
\end{align*}
Then by a union bound over all vertices in $V_2$ and all choices $\{i,j\}$ with $1\le i<j\le r$, we have $\delta^\ast(B_{\pi})\ge \left(1-\frac{1}{r}+\frac{\gamma}{2}\right)n$ with probability at least $1-2n{r\choose 2}\exp\{-(\gamma/32)^2n\}$.
\end{proof}

\end{proof}



\bibliographystyle{abbrv}
\bibliography{ref}

\begin{appendices}
\section{Proof of Lemma \ref{lem: Pham}}\label{appendix}
\begin{defn}
For a finite set $X$ and its subset $S$, let $\langle S\rangle=\{Y: S\subseteq Y\subseteq X\}$ and $\mathcal{F}\subseteq 2^X$. We say that a probability measure $\rho$ on $\mathcal{F}$ is \emph{$q$-spread} if
$$\rho(\mathcal{F}\cap \langle S\rangle)\le q^{|S|}~\text{for any}~S\subseteq X.$$
\end{defn}

\begin{theorem}[\cite{Pham2}, Theorem 1.9]\label{lem:Pham2}
Let $r\ge 2$ and $d>0$. There exists a constant $\varepsilon>0$ such that the following holds for any sufficiently large $n$. Let $G=(V_1,V_2,\ldots, V_r,E)$ be a balanced $r$-partite graph with each vertex class of size $n$. Suppose that $(V_i,V_j)$ is $(\varepsilon,d_{ij})$-super regular with $d_{ij}\ge d$ for all $i\neq j$. Let $H$ be the $r$-clique complex of $G$ and $\mathcal{K}\subseteq 2^{H}$ be the set of $K_r$-factors in $G$. Then there is a $O(1/n^{r-1})$-spread probability measure on $\mathcal{K}$.
\end{theorem}

The following theorem can be deduced from Theorem 1.6 of \cite{Frankston}.

\begin{theorem}[\cite{Frankston}]\label{lem:Frankston}
There exists a constant $C>0$ such that the following holds. Let $X$ be a finite set and $\mathcal{F}\subseteq 2^X$. If there is a $q$-spread probability measure on $\mathcal{F}$, then whp $X(\min\{Cq\log |X|,1\})$ has a subset which is an element of $\mathcal{F}$.
\end{theorem}

\begin{proof}[Proof of Lemma \ref{lem: Pham}]
Given an integer $r\ge 3$ and $d>0$, we choose $$\frac{1}{n}\ll\frac{1}{C}\ll\varepsilon\ll d, \frac{1}{r}.$$
Recall that $G=(V_1,V_2,\ldots, V_r,E)$ is a balanced $r$-partite graph with each vertex class of size $n$ and $(V_i,V_j)$ is $(\varepsilon,d_{ij})$-super regular with $d_{ij}\ge d$ for all $i\neq j$.
Let $H$ be the $r$-clique complex of $G$ and $\mathcal{K}\subseteq 2^{H}$ be the set of $K_r$-factors in $G$.
By Theorem \ref{lem:Pham2}, there is a $O(1/n^{r-1})$-spread probability measure on $\mathcal{K}$.
Let $p_0=\Omega(\log n/n^{r-1})$.
Then by applying Theorem \ref{lem:Frankston} with $X=H$, $\mathcal{F}=\mathcal{K}$ and $q=O(1/n^{r-1})$, we obtain that whp $H(p_0)$ has a subset which is an element of $\mathcal{K}$, that is, a $K_r$-factor in $G$.
Let $p\ge Cn^{-2/r}(\log n)^{1/{r\choose 2}}$.
By Lemma \ref{lem: clique complex},  whp $G(p)$ contains a $K_r$-factor.
\end{proof}

\end{appendices}

\end{document}